\documentclass{article}

\usepackage[english]{babel}
\usepackage{multirow}
\usepackage{blindtext}
\usepackage{booktabs}
\usepackage{siunitx}
\usepackage{graphicx}
\usepackage{subcaption}
\usepackage{float} 
\usepackage{enumitem} 

\usepackage{natbib}
\usepackage{amssymb}

\usepackage[letterpaper,top=2cm,bottom=2cm,left=3cm,right=3cm,marginparwidth=1.75cm]{geometry}

\usepackage{amsmath}
\usepackage{graphicx}
\usepackage[colorlinks=true, allcolors=blue]{hyperref}

\begin{document}
\title{PyClustrPath: An efficient Python package for generating clustering paths with GPU acceleration\footnotemark[1]}
\author{Hongfei Wu\footnotemark[2], \quad Yancheng Yuan \footnotemark[3]}
\maketitle
\renewcommand{\thefootnote}{\fnsymbol{footnote}}
\footnotetext[1]{Corresponding author: Yancheng Yuan.\\
\textbf{Funding}: The work of Yancheng Yuan
was supported by the Hong Kong Research Grants Council (Project No. 25305424) and the Research
Center for Intelligent Operations Research at The Hong Kong Polytechnic University.}
\footnotetext[2]{Department of Data Science and Artificial Intelligence, The Hong Kong Polytechnic University, Hung Hom, Hong Kong ({\tt hung-fei.wu@connect.polyu.hk}).}
\footnotetext[3]{Department of Data Science and Artificial Intelligence, The Hong Kong Polytechnic University, Hung Hom, Hong Kong ({\tt yancheng.yuan@polyu.edu.hk}).}
\renewcommand{\thefootnote}{\arabic{footnote}}

\begin{abstract}
Convex clustering is a popular clustering model without requiring the number of clusters as prior knowledge. It can generate a clustering path by continuously solving the model with a sequence of regularization parameter values. This paper introduces {\it PyClustrPath}, a highly efficient Python package for solving the convex clustering model with GPU acceleration. {\it PyClustrPath} implements popular first-order and second-order algorithms with a clean modular design. Such a design makes {\it PyClustrPath} more scalable to incorporate new algorithms for solving the convex clustering model in the future.  We extensively test the numerical performance of {\it PyClustrPath} on popular clustering datasets, demonstrating its superior performance compared to the existing solvers for generating the clustering path based on the convex clustering model. The implementation of {\it PyClustrPath} can be found at: \href{https://github.com/D3IntOpt/PyClustrPath}{https://github.com/D3IntOpt/PyClustrPath}. 
\end{abstract}


\section{Introduction}
Convex clustering model \citep{pelckmans2005convex,lindsten2011clustering,hocking2011clusterpath} has become popular due to its superior and robust performance in clustering unstructured data, especially when the prior knowledge of the number of clusters is absent. For a given collection of $n$ data points with $d$ features $A=\{\mathbf{a}_1, \mathbf{a}_2, \ldots, \mathbf{a}_n\} \subseteq \mathbb{R}^{d}$, the convex clustering model can generate a clustering path by continuously solving the following strongly convex model for a sequence of values of $\gamma$:
\begin{equation}
\tag{CCM}
\label{model: CCM}
\min _{x_1, \dots, x_n \in \mathbb{R}^{d}} ~ \frac{1}{2} \sum_{i=1}^n\left\|\mathbf{x}_i-\mathbf{a}_i\right\|_2^2+\gamma \sum_{i<j} w_{i j}\left\|\mathbf{x}_i-\mathbf{x}_j\right\|_q,
\end{equation}
where $w_{i j}=w_{j i} \geq 0$ are given weights depending on the input data $A$, $\gamma>0$ is a tuning parameter that controls the strength of the fusion penalty, and $\|\cdot\|_q$ is the vector $q$-norm ($q \geq 1$). For any fixed $\bar{\gamma} \geq 0$, after obtaining the unique solution $\mathbf{x}_i^{*}(\bar{\gamma}) ~ (1 \leq i \leq n)$ to the model (\ref{model: CCM}), we will assign points $\mathbf{a}_i$ and $\mathbf{a}_j$ into the same cluster if $\mathbf{x}_i^{*}(\bar{\gamma})=\mathbf{x}_j^{*}(\bar{\gamma})$. When $\gamma = 0$, we have $\mathbf{x}_i^*(0) = \mathbf{a}_i$, thus the convex clustering model gives $n$ clusters if the input data points are distinct. As we increase the value of $\gamma$, some $\mathbf{x}^*_i(\gamma)$ will become identical due to the fusion penalty terms in the model. If the undirected weighted graph specified by $\{w_{ij}\}_{1 \leq i < j \leq n}$ is connected, it holds that $\mathbf{x}^*_i(\hat{\gamma}) = \mathbf{x}^*_j(\hat{\gamma}) ~ \forall 1 \leq i \leq j \leq n$ for sufficiently large $\hat{\gamma} > 0$. Moreover, \citet{chi2015splitting} have proved that $\mathbf{x}^*(\gamma)$ is a continuous function of $\gamma$ for $\gamma > 0$. In practice, we usually solve \eqref{model: CCM} for a sequence of values for the parameter $\gamma$, i.e., $0\leq \gamma_1 < \gamma_2 < \cdots < \gamma_T < +\infty$, and obtain a clustering path of the data points. 

Motivated by the superior empirical performance of \eqref{model: CCM}, extensive investigations on the theoretical and algorithmic perspectives of the convex clustering model have been done. We briefly mention some advances in the convex clustering model below, a more comprehensive survey of the convex clustering model can be found in \citep{feng2023review} and the references therein. On the one hand, the perfect (and the nontrivial) cluster recovery guarantees have been established for \eqref{model: CCM} under some practical assumptions \citep{zhu2014convex,panahi2017clustering,sun2021convex,jiang2020recovery,chiquet2017fast,chi2019recovering,dunlap2021sum}. The statistical properties of the convex clustering model have also been extensively investigated \citep{tan2015statistical, Radchenko2017,chi2020provable}. On the other hand, several optimization algorithms have been designed for solving the convex clustering model. In particular, \cite{chi2015splitting} proposed to solve \eqref{model: CCM} using the alternating direction method of multipliers (ADMM) and fast alternating minimization algorithm (AMA). They also developed a popular R package {\it cvxclustr}\footnote{https://github.com/echi/cvxclustr} based on these two algorithms.  Later, \cite{yuan2018efficient} proposed a semismooth Newton based augmented Lagrangian method (SSNAL), which achieved superior numerical performance for solving large-scale convex clustering model \eqref{model: CCM}, especially when solutions of high accuracy are required. They implemented the SSNAL algorithm and released a Matlab package {\it ConvexClustering}\footnote{https://blog.nus.edu.sg/mattohkc/softwares/convexclustering/}.  Currently, neither of the two popular packages supports computation on a Graphics Processing Unit (GPU).

\begin{table}[ht]
\centering
\caption{Feature comparison among the popular packages for solving the convex clustering model.}
\label{table: package-comp}
\begin{tabular}{@{}lcccc @{}} 
\toprule
       & AMA       & ADMM  & SSNAL   & GPU Support  \\ \midrule
\textbf{\it cvxclustr} ~ \citep{chi2015splitting}& $\checkmark$ & $\checkmark$ & $\times$ & $\times$ \\ 
\textbf{\it ConvexClustering} ~ \citep{sun2021convex}      & $\checkmark$ & $\checkmark$ & $\checkmark$ & $\times$ \\ 
\textbf{\it PyClustrPath} (\textbf{This paper}) & $\checkmark$ & $\checkmark$ & $\checkmark$ & $\checkmark$\\ 
\bottomrule
\end{tabular}
\end{table}

This paper introduces a new Python package {\it PyClustrPath} for efficiently generating clustering paths by solving \eqref{model: CCM}, which supports both CPU and GPU computation simultaneously. We summarize some key features among the three packages: {\it cvxclustr}, {\it ConvexClustering}, and {\it PyClustrPath} in Table \ref{table: package-comp}. Note that the additions, matrix-vector multiplications, and matrix factorizations dominate the computational cost for all three popular algorithms: ADMM, fast AMA, and SSNAL. Therefore, we can take advantage of GPU for accelerating these algorithms. We summarize the main contributions of {\it PyClustrPath} as follows:

\begin{itemize}
    \item We design and develop a new Python package for efficiently generating clustering paths by solving \eqref{model: CCM}, which supports computation on both CPU and GPU simultaneously. 

    \item We implement three popular and efficient algorithms in {\it PyClustrPath}: ADMM, fast AMA, and SSNAL. 

    \item We design the package {\it PyClustrPath} in a well-structured modular manner, making it flexible and scalable to incorporate new algorithms for solving \eqref{model: CCM} in the future. 

    \item We develop a visualization module in {\it PyClustrPath}, providing a user-friendly way to visualize the generated clustering path by \eqref{model: CCM}.

    \item We extensively test {\it PyClustrPath} 's numerical performance for generating clustering paths for several popular benchmark datasets. The results demonstrate the efficiency and robustness of {\it PyClustrPath}.
\end{itemize}

The rest of the paper is organized as follows. Section 2 describes the architecture and design of the {\it PyClustrPath} package. Section 3 includes instructions and demos for using the {\it PyClustrPath} package. Detailed results of the numerical performance of the {\it PyClustrPath} package for solving the convex clustering model will be shown in Section 4. We conclude the paper and discuss some future research directions in Section 5.  

\section{Package architecture and design}
The architecture of \textit{PyClustrPath} is designed to ensure efficiency, modularity, scalability, and user-friendliness. The package aims to simplify the development, integration, and application of convex clustering algorithms while ensuring high computational efficiency. This section describes its dependencies, workflow, and key features.

\subsection{Dependencies}

To ensure high performance and flexibility, \textit{PyClustrPath} builds on a set of well-established frameworks and libraries. These dependencies are chosen not only for their computational efficiency but also for their compatibility with modern hardware platforms.

\textit{PyClustrPath} is built for Python 3.10+ and leverages the PyTorch framework~\citep{paszke2019pytorch} to support computations on both GPU and CPU platforms. PyTorch is chosen due to its widespread adoption in machine learning and scientific computing communities, its ease of integration with CUDA for GPU acceleration, and its robust support for tensor operations. These features are crucial for handling high-dimensional and large-scale datasets in clustering tasks. Additionally, \textit{PyClustrPath} adopts \textit{cholespy}~\citep{Nicolet2021Large}, a specialized library for efficient sparse Cholesky factorizations.
Together, these dependencies ensure {\it PyClustrPath} to solve large-scale convex clustering problems on diverse hardware platforms efficiently.

\subsection{Architecture}

\begin{figure}[H]
    \centering
    \includegraphics[width=0.8\linewidth]{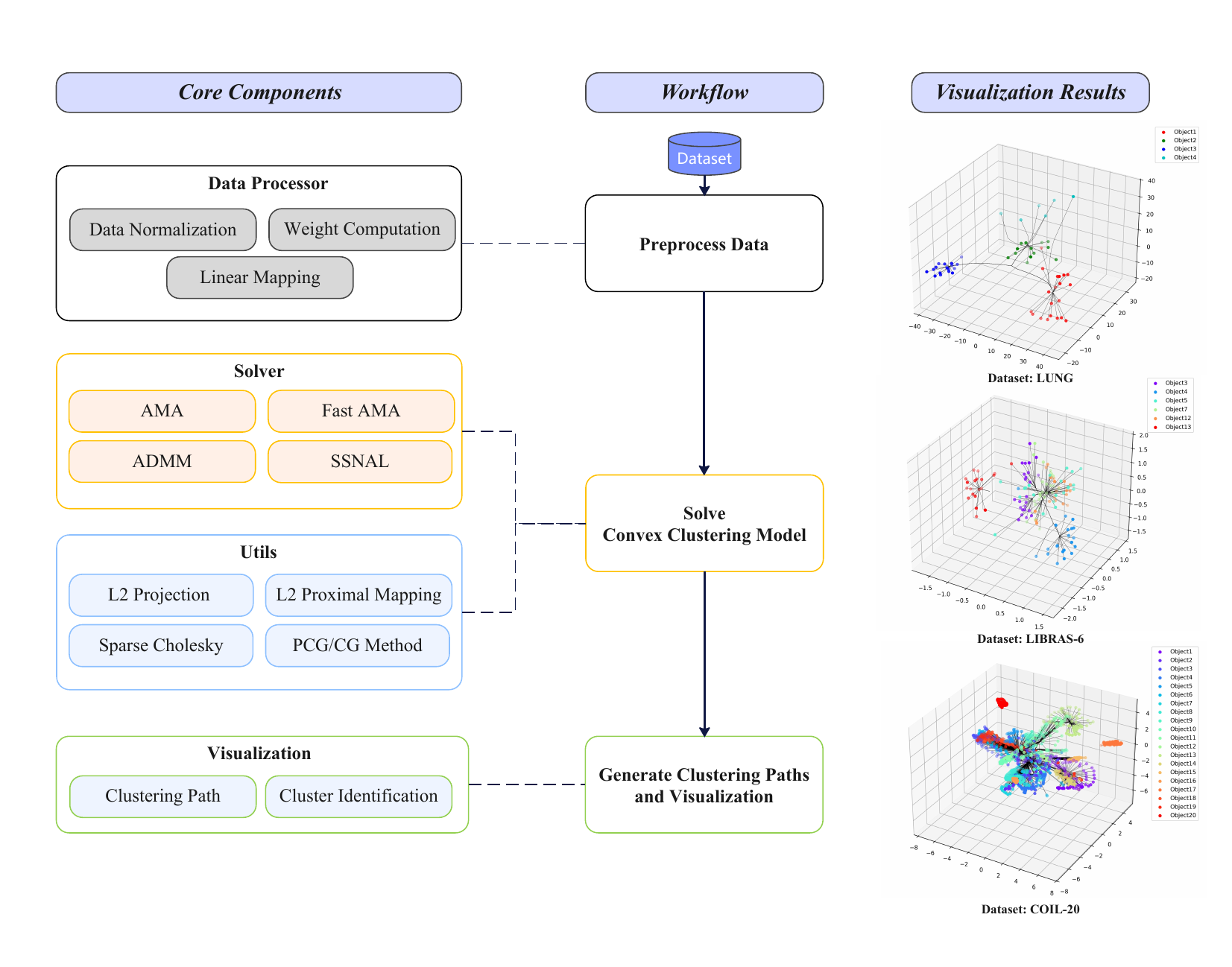}
    \caption{Components, Workflow and Visualization results.}
    \label{fig:Components-Workflow}
\end{figure}

The architecture of \textit{PyClustrPath} is organized into four core components: Data Processor, Solver, Utilities, and Visualization, as illustrated in Figure~\ref{fig:Components-Workflow}. Each component is designed to perform specific tasks in the workflow:
(1) \textbf{Data Processor} handles data preprocessing, including the computation of weights, the construction of linear maps required for solving \eqref{model: CCM}, and so on.
(2) \textbf{Solver} serves as the computational core of the package, offering users a selection of popular algorithms. Users can freely choose the algorithm according to their preferences and needs.
(3) \textbf{Utils} provides efficient computation subroutines which are required in the optimization algorithms implemented in the \textbf{Solver} component.
(4) \textbf{Visualization} provides a convenient interface for displaying the clustering paths.

\subsection{Features}

The design of \textit{PyClustrPath} incorporates several key features that distinguish it from existing packages for solving the convex clustering model:

\begin{itemize}
    \item \textbf{Comprehensive algorithm selection.} \textit{PyClustrPath} supports multiple popular and efficient convex clustering algorithms, offering users a broad range of choices to address diverse clustering tasks. It integrates {\it cholespy} package for sparse Cholesky decomposition and implements an efficient preconditioned conjugate gradient (PCG)~\citep{nocedal2006numerical} to ensure efficiency and scalability of {\it PyClustrPath} for solving clustering tasks with huge-scale datasets.

    \item \textbf{Efficient computation with GPU supports.} By leveraging PyTorch's supports for the GPU computation, combined with custom CUDA extensions, \textit{PyClustrPath} achieves significant speedups over existing CPU-based implementations. A detailed comparison of computational performance is provided in Section~\ref{sec:experiments}.

    \item \textbf{Modular design.} Inspired by the principles of scikit-learn~\citep{buitinck2013api}, the package adopts a modular structure where all convex clustering algorithms inherit from a unified base class and follow a consistent API interface. This modularity enables users to easily switch between algorithms and gives researchers the flexibility to extend the package with new algorithms.
    
\end{itemize}

The combination of these features ensures that \textit{PyClustrPath} is not only a powerful package to solve \eqref{model: CCM} but also a versatile platform to advance research in related fields.

\section{Instructions and demos for using {\it PyClustrPath}}

\textit{PyClustrPath} allows users to freely select the optimization algorithm and customize various parameters in the solver, such as the maximum number of iterations, tolerance levels, whether to utilize GPU acceleration, and so on. This section describes the key steps for using the package by providing several demo examples of different clustering datasets. Detailed documentation is available on our GitHub repository for the {\it PyClustrPath} package\footnote{ \url{https://github.com/D3IntOpt/PyClustrPath}}.

We demo the usage of {\it PyClustrPath} on several real-world datasets: LIBRAS \citep{misc_libras_movement_181}, LIBRAS-6 \citep{misc_libras_movement_181}, COIL-20 \citep{nene1996columbia}, LUNG \citep{lee2010biclustering}, and MNIST \citep{lecun1998gradient}. These datasets cover diverse domains, including time series, image, and gene expression data, which highlight the flexibility and robustness of \textit{PyClustrPath} for solving the convex clustering model \eqref{model: CCM}. Table~\ref{table:dataset} summarizes the key characteristics of the datasets.

\begin{table}[ht]
\centering
\caption{Key characteristics of the real datasets. Here, \(d\) represents the number of features, \(n\) is the number of samples, \(K\) is the ground-truth number of clusters, and \(\gamma\) denotes the range of the regularization parameter used for generating clustering paths by the convex clustering model \eqref{model: CCM}.}
\label{table:dataset}
\begin{tabular}{@{}lccccc@{}} 
\toprule
Dataset      & Type       & \(d\)  & \(n\)   & \(K\) & $\gamma$ \\ \midrule
LIBRAS-6     & time series & 90   & 144    & 6  & [0.45 : 0.09] \\
LIBRAS       & time series & 90   & 360    & 15 & [0.60 : 0.30] \\
MNIST        & image       & 784  & 10000  & 10 & [0.69 : 0.48] \\
COIL-20      & image       & 1024 & 1440   & 20 &[7.10 : 0.64] \\
LUNG         & gene        & 12625& 56     & 4  &[0.93 : 0.057] \\ \bottomrule
\end{tabular}
\end{table}

Next, we demo the usage of {\it PyClustrPath} for generating the clustering path on these datasets.

\subsection{A demo for the LIBRAS-6 dataset}
\begin{figure}[H]
    \centering
    \includegraphics[width=\linewidth]{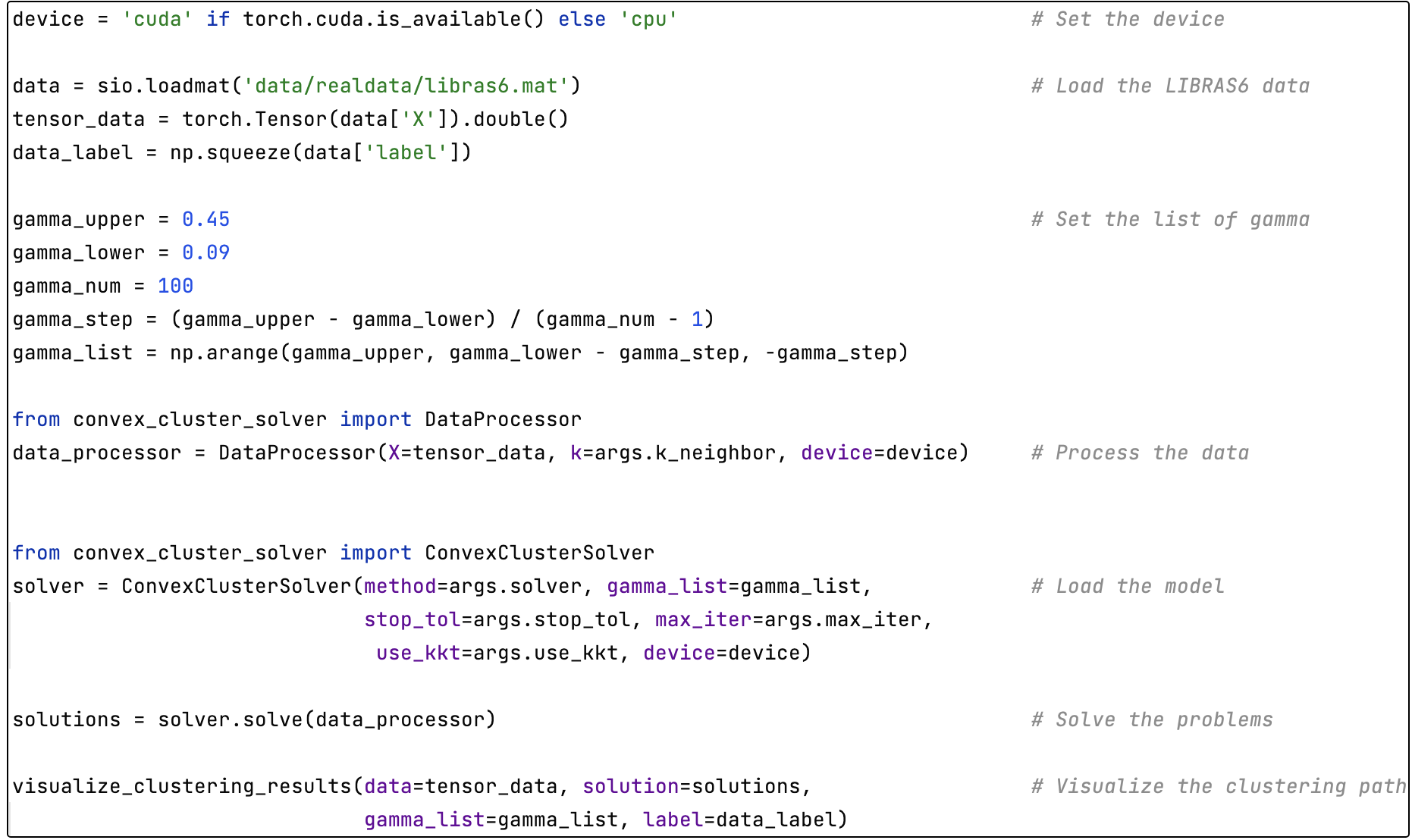}
    \caption{Demo codes for LIBRAS-6 dataset.}
    \label{fig:code-libras6}
\end{figure}

\begin{figure}[H]
    \centering
    \includegraphics[width=0.7\linewidth]{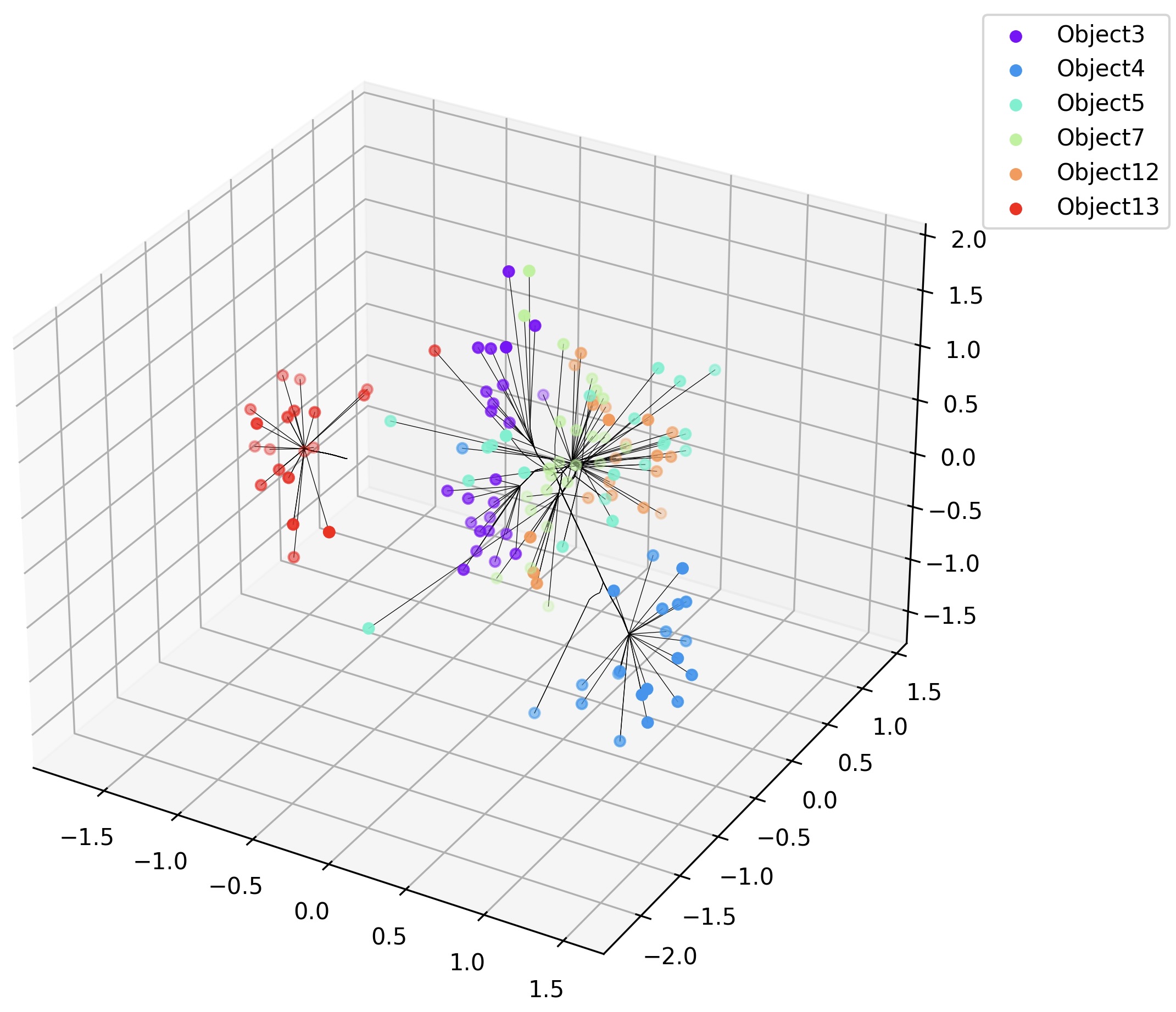}
    \caption{A generated clustering path for the LIBRAS-6 dataset.}
    \label{fig:visual_libras6}
\end{figure}

Figure \ref{fig:code-libras6} provides the demo codes for generating clustering path by {\it PyClustrPath} on the LIBRAS-6 dataset, including necessary annotations for explaining the key steps. The generated clustering path is shown in Figure \ref{fig:visual_libras6}.

\subsection{A demo for the COIL-20 dataset}
\begin{figure}[H]
    \centering
    \includegraphics[width=\linewidth]{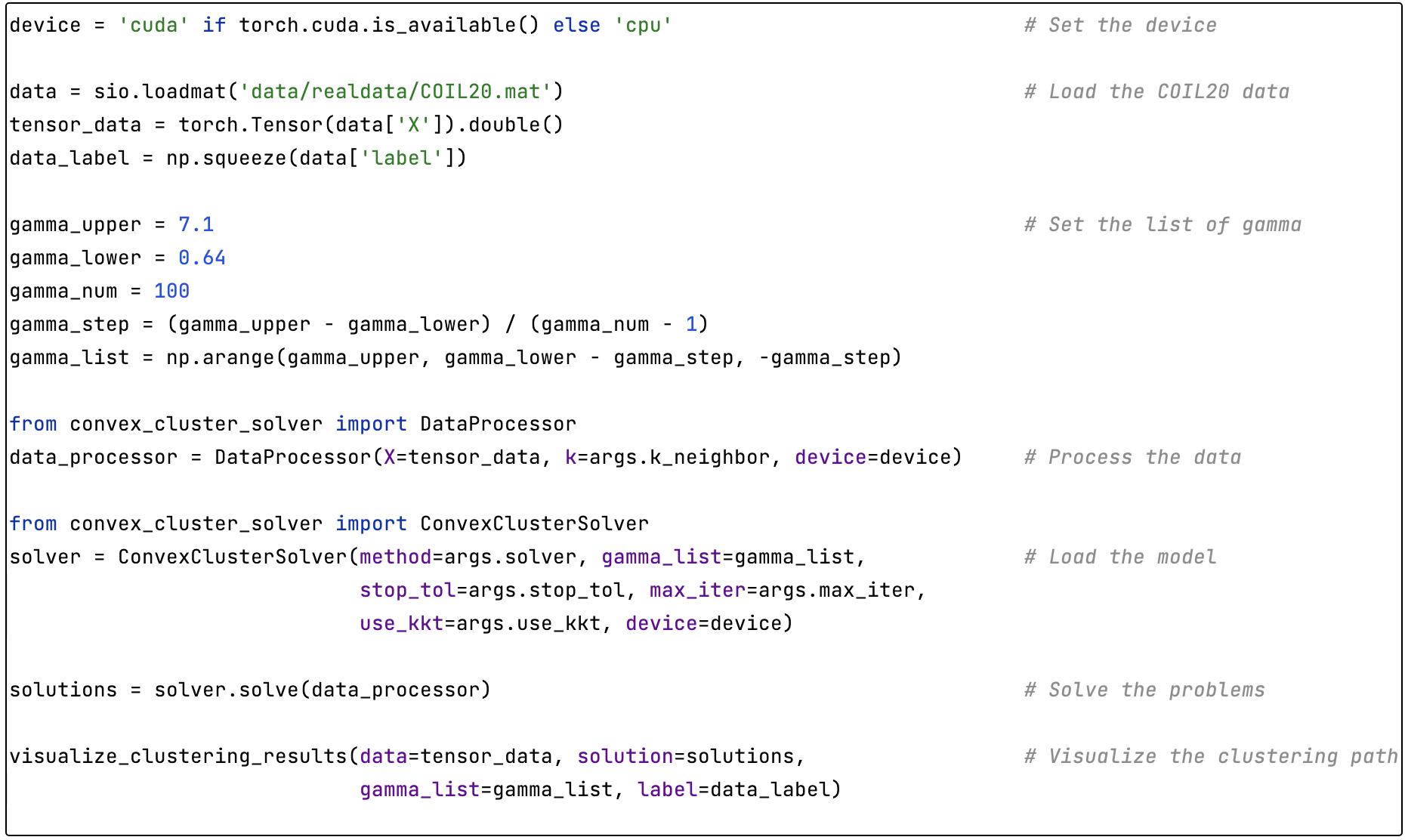}
    \caption{Demo codes for the COIL-20 dataset.}
    \label{fig:code-coil20}
\end{figure}

\begin{figure}[H]
    \centering
    \includegraphics[width=0.7\linewidth]{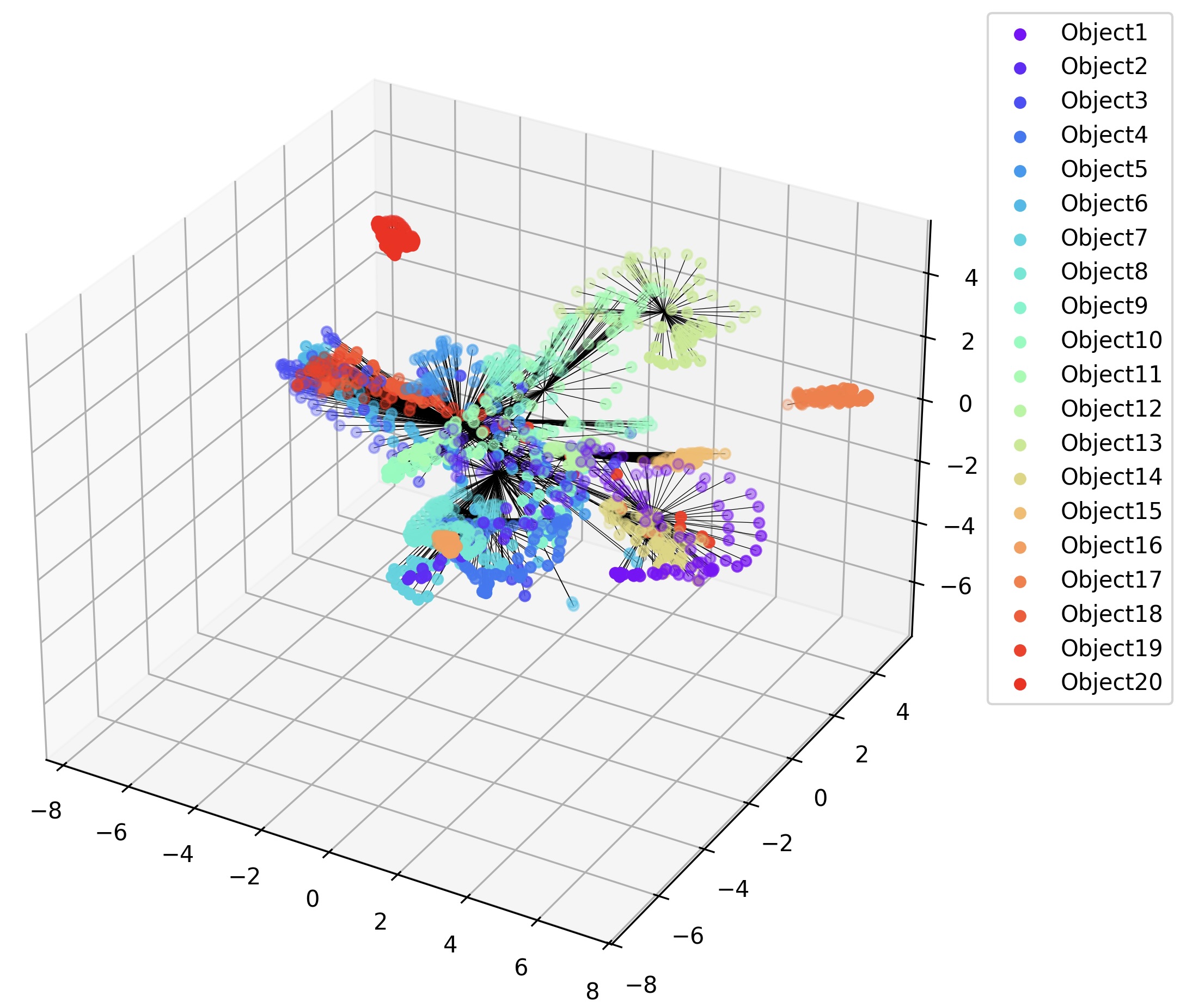}
    \caption{A generated clustering path for the COIL-20 dataset.}
    \label{fig:visual_coil20}
\end{figure}

Figure \ref{fig:code-coil20} provides the demo codes for generating clustering path by {\it PyClustrPath} on the COIL-20 dataset. A visualization of the generated clustering path is in Figure \ref{fig:visual_coil20}.

\subsection{A demo for the LUNG dataset}
\begin{figure}[H]
    \centering
    \includegraphics[width=\linewidth]{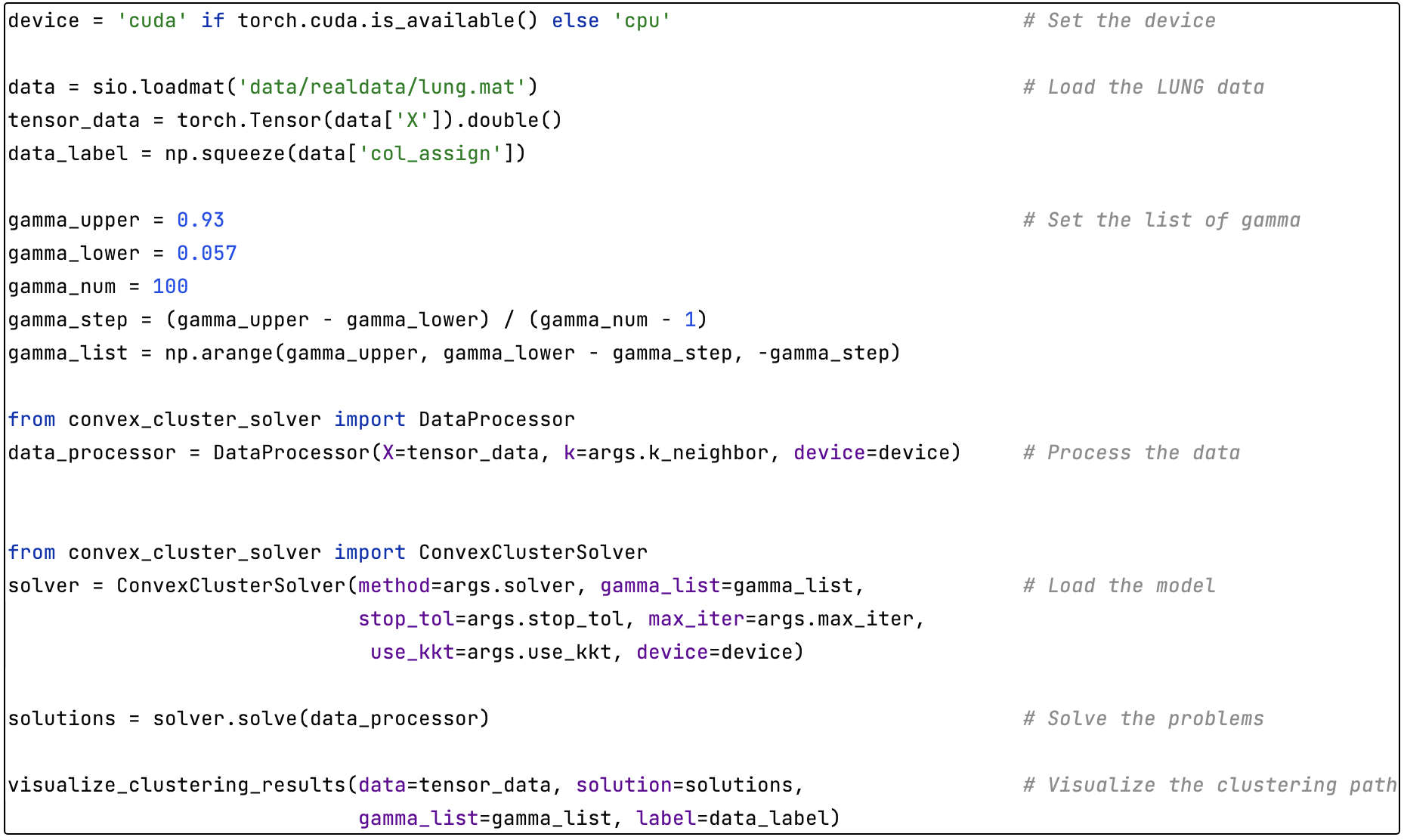}
    \caption{Demo codes for the LUNG dataset.}
    \label{fig:code-lung}
\end{figure}

\begin{figure}[H]
    \centering
    \includegraphics[width=0.7\linewidth]{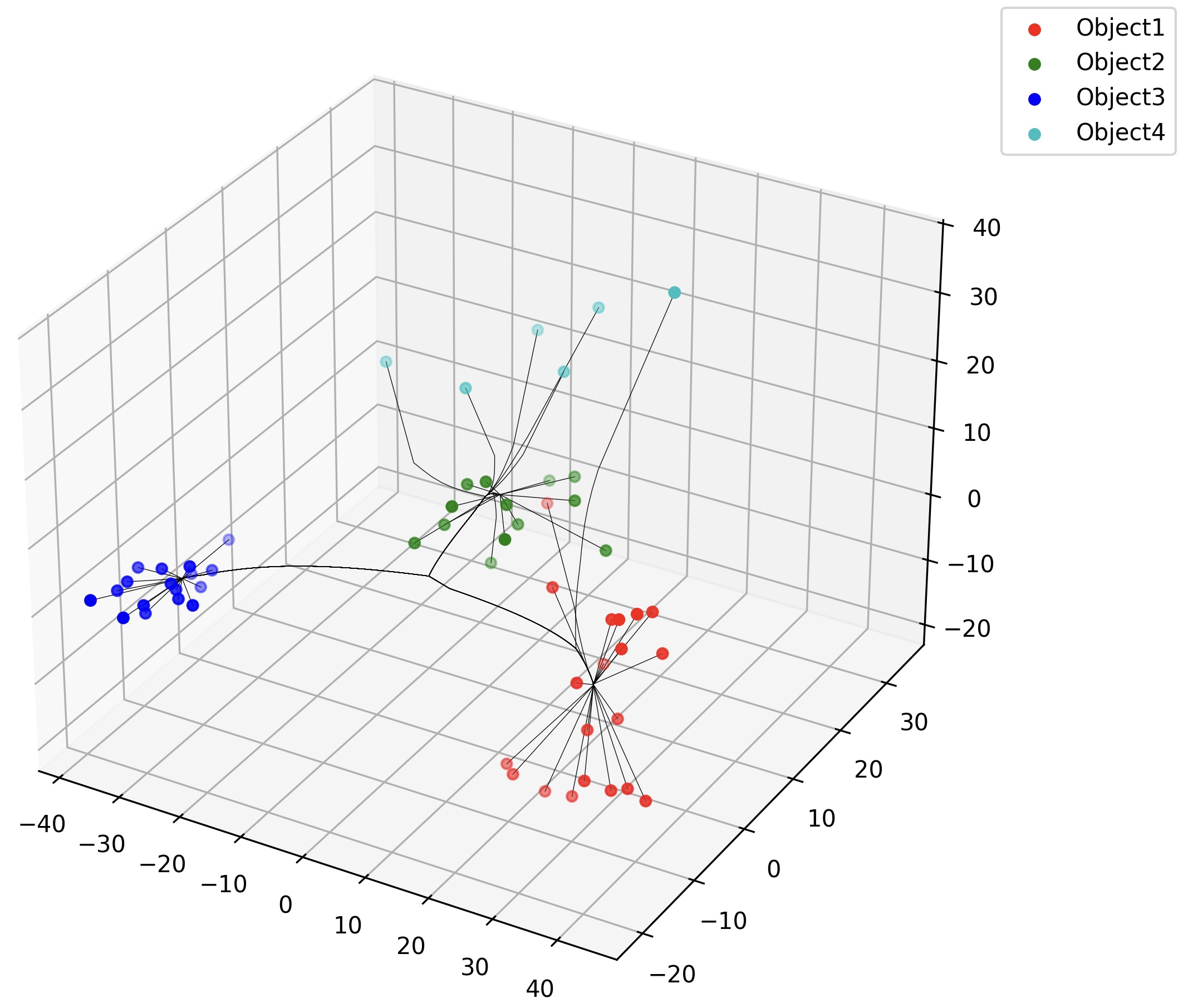}
    \caption{A generated clustering path for the LUNG dataset.}
    \label{fig:visual_lung}
\end{figure}

Figure \ref{fig:code-lung} provides the demo codes for generating clustering path by {\it PyClustrPath} on the LUNG dataset. The generated clustering path can be found in Figure \ref{fig:visual_lung}.

\subsection{A demo for the MNIST dataset}
\begin{figure}[H]
    \centering
    \includegraphics[width=\linewidth]{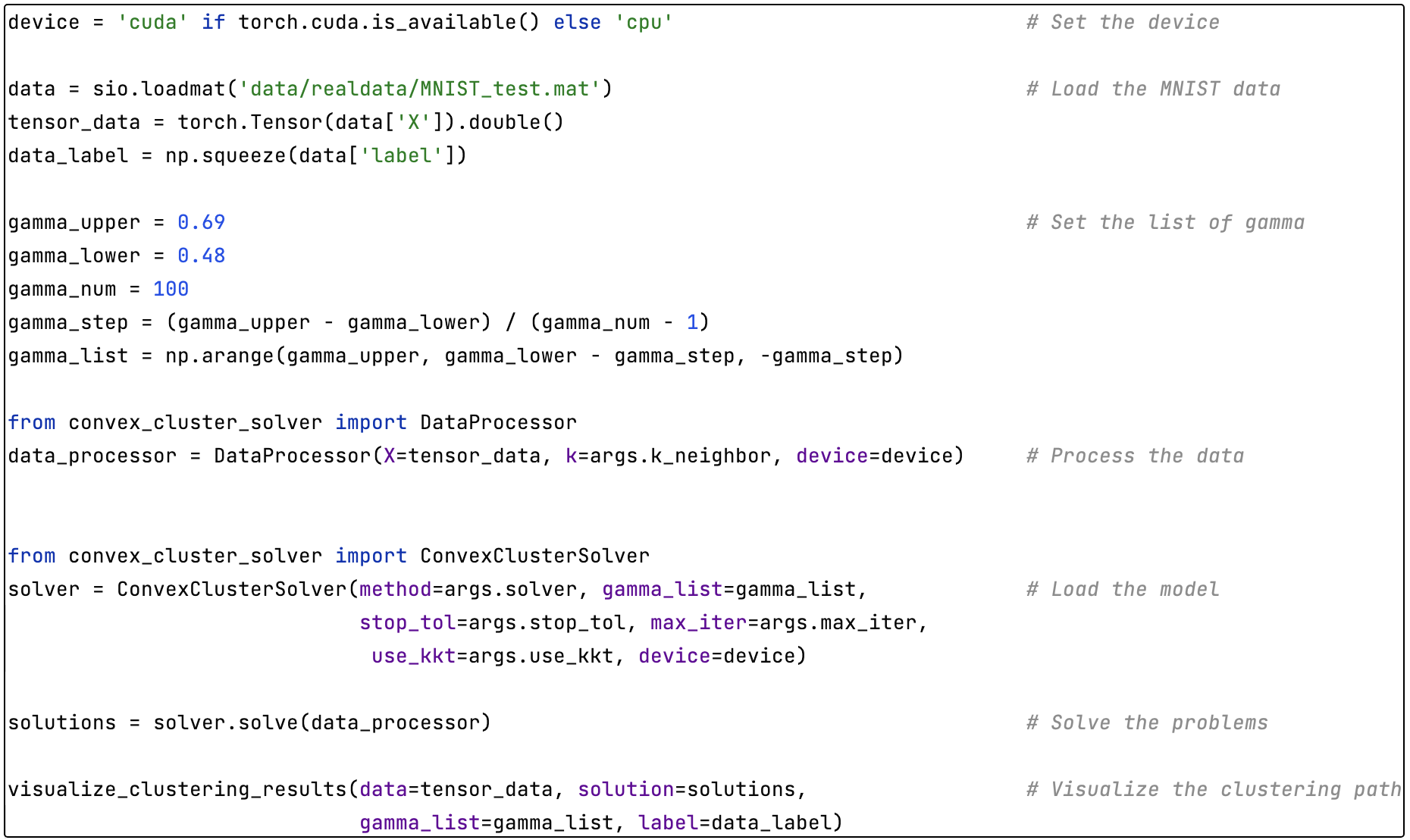}
    \caption{Usage Instructions for MNIST dataset.}
    \label{fig:code-mnist}
\end{figure}

\begin{figure}[H]
    \centering
    \includegraphics[width=0.7\linewidth]{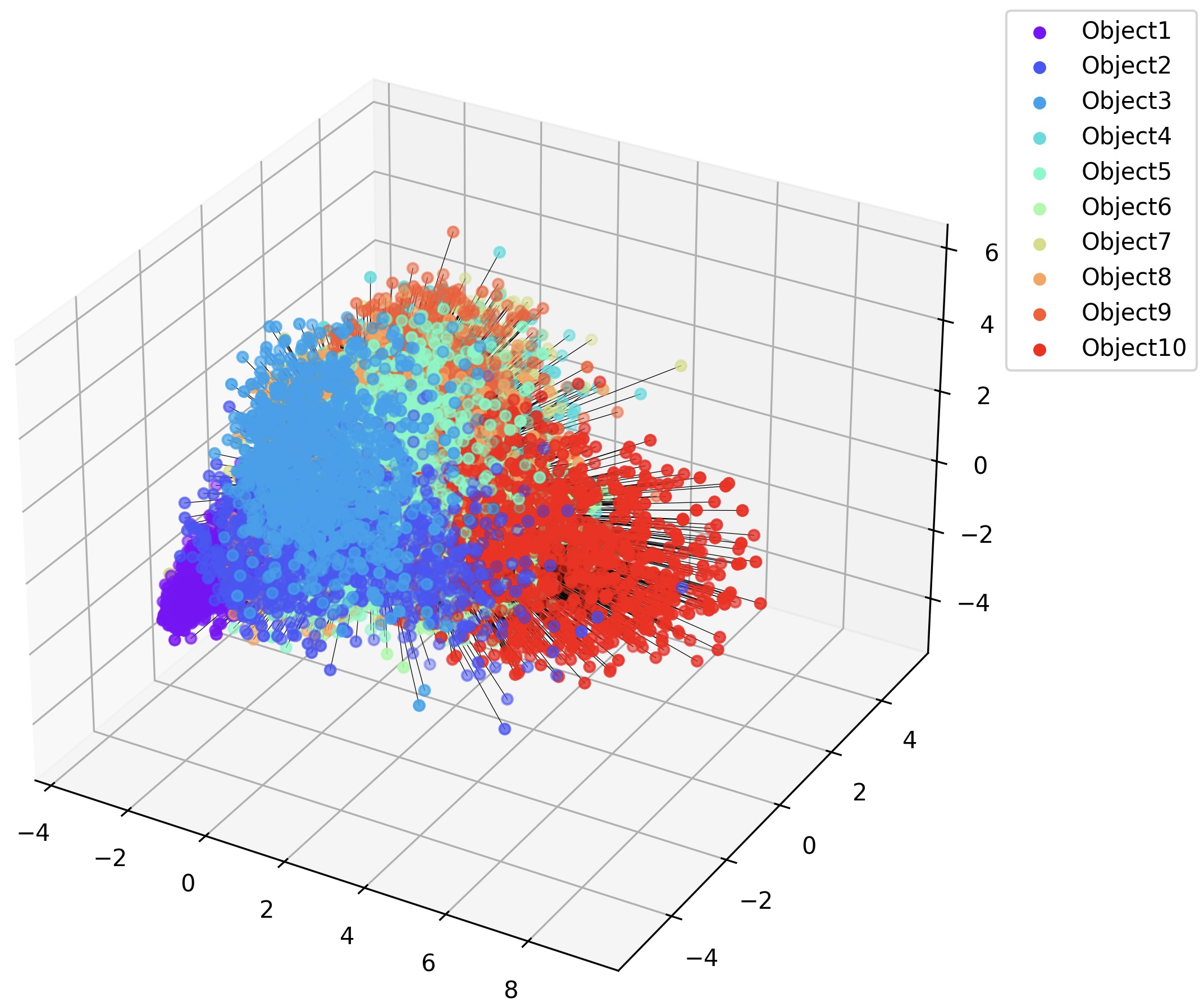}
    \caption{A generated clustering path for the MNIST dataset.}
    \label{fig:visual_mnist}
\end{figure}

Figure \ref{fig:code-mnist} provides the demo codes for generating clustering path by {\it PyClustrPath} on the MNIST dataset. The corresponding generated clustering path can be found in Figure \ref{fig:visual_mnist}.

\section{Numerical experiments} \label{sec:experiments}
This section provides comprehensive evaluations of  {\it PyClustrPath} and the two popular packages {\it cvxclustr} ~ \citep{chi2015splitting} and {\it ConvexClustering} ~ \citep{yuan2018efficient,sun2021convex} for solving \eqref{model: CCM} on five commonly used benchmark datasets: LIBRAS-6, LIBRAS, COIL-20, LUNG, and MNIST. All tests were conducted on a platform with an Intel(R) Xeon(R) Platinum 8480C CPU and a single NVIDIA GeForce RTX 4090 Graphics card. The comparisons focus on computational efficiency across algorithms and platform settings (CPU and GPU).

To evaluate the performance, we generate clustering paths for each dataset by solving problem \eqref{model: CCM} for a sequence of 100 values of $\gamma > 0$. This corresponds to solving 100 convex clustering problems per dataset. We follow \citep{yuan2021dimension} to use the relative duality gap 
\[
\eta = \frac{|f_{p} - f_d|}{1 + |f_p| + |f_d|} \leq \epsilon
\]
to terminate all the algorithms, where $f_p$ and $f_d$ are the objective values for the primal problem \eqref{model: CCM} and the corresponding dual problem, respectively. 

To ensure reliability and consistency, we employ the following experimental strategies:
(1) \textbf{Baseline time setting}: The baseline runtime \( T \) is typically determined by the best-performing algorithm on the GPU due to its parallel computing efficiency. To compare all methods within a reasonable timeframe while still capturing meaningful performance differences, we set the maximum runtime for each algorithm to \( 10T \) on the corresponding dataset.
(2) \textbf{Annotation of results}: For cases where 100 problems could not be solved within the running time of \( 10T \), we annotate the results with \textit{*(s)} to indicate a number of \textit{s} problems successfully solved within this time limit. For example, on the COIL-20 dataset, the fast AMA algorithm in {\it PyClustrPath} (GPU) solves only 14 problems within 482.85 seconds, which is denoted as \( 482.85^*(14) \). We summarize these numerical results in Table~\ref{tab:computation_time_1e-6}, which demonstrates the superior performance of {\it{PyClustrPath}} with GPU acceleration.

\begin{table}[htbp]
\centering
\caption{Computation time (seconds) for various datasets with tolerance $\epsilon = 10^{-6}$.  
``$\times$'' means that the algorithm or platform is unsupported or failed to process the dataset. Results annotated with \textit{*(s)} indicate a number of \textit{s} problems solved within the set runtime threshold \( 10T \).}
\label{tab:computation_time_1e-6}
\setlength{\tabcolsep}{14pt}
\renewcommand{\arraystretch}{1.2} 
\resizebox{\textwidth}{!}{%
    \begin{tabular}{
      l|
      l|
      S[table-format=4.2]
      S[table-format=4.2]
      S[table-format=4.2]
      S[table-format=4.2]
    }
    \toprule
    {Dataset} & {Method} & {{\it cvxclustr} (CPU)} & {{\it ConvexClustering} (CPU)} &{{\it PyClustrPath} (CPU)} &{{\it PyClustrPath} (GPU)}\\
    \midrule
    \multirow{3}{*}{LIBRAS-6} 
                & SSNAL & $\times$  & 40.42  & 22.73 & \textbf{8.94} \\
               & ADMM & 96.06*(15) & 98.64*(92)& 92.48 & 59.35 \\
               & Fast AMA & 64.98 & 90.06*(55) & 89.95*(36) & 83.99 \\
    \midrule
    \multirow{3}{*}{LIBRAS}   
                & SSNAL & $\times$  & 39.91 & 17.82 & \textbf{5.70} \\
               & ADMM & 197.49*(1) & 42.25 & 20.04 & 7.06 \\
               & Fast AMA & 62.51*(12) & 59.35*(4) & 60.72*(5) & 57.45*(36) \\
    \midrule
    \multirow{3}{*}{COIL-20}
                & SSNAL & $\times$ & 469.28*(48) & 472.14*(63)  & \textbf{46.78} \\
               & ADMM & $\times$ & 472.28*(47) & 469.67*(57) & 102.90 \\
               & Fast AMA & 10868.27*(0) & 5278.81*(1) & 5804.21*(1) & 482.85*(14) \\
    \midrule
    \multirow{3}{*}{LUNG}       
                & SSNAL & $\times$ & 690.50*(89) & 660.14 & 67.34 \\
               & ADMM & 733.45*(18) & 621.56*(85) & 735.05*(86) & 665.87*(96) \\
               & Fast AMA & 676.45*(58) & 1320.53*(4) & 1287.75*(4) & \textbf{65.08} \\
    \midrule
    \multirow{3}{*}{MNIST}
                & SSNAL & $\times$ & 8920.74*(17) & 8928.26*(95) & \textbf{889.87} \\
               & ADMM & $\times$ & 8909.42*(7) & 8965.12*(84) & 915.35 \\
               & Fast AMA & 41606.83*(0) & 14869.85*(1) & 16 601.79*(1)  & 3278.37 \\
    \bottomrule
    \end{tabular}
}
\end{table}

\subsection{Performance across datasets}

The efficiency of the GPU acceleration is evident across datasets of varying scales, including small datasets such as LIBRAS-6 and LIBRAS, as well as large and more computationally demanding datasets like COIL-20, LUNG, and MNIST.

\begin{itemize}
    \item \textbf{LIBRAS-6 and LIBRAS}. For the small LIBRAS-6 dataset, GPU acceleration provides noticeable improvements.  SSNAL with the GPU acceleration outperforms the best CPU version (achieved by {\it PyClustrPath}), achieving a nearly 2.5x speedup.
    On the LIBRAS dataset, which is slightly larger, the GPU implementation maintains superior performance. SSNAL on GPU achieves a 3.1x speedup. Similarly, ADMM and fast AMA benefit significantly from the GPU acceleration.

    \item \textbf{COIL-20 and LUNG}. 
    As the scale of the dataset increases, the advantage of GPU becomes more pronounced. SSNAL with the GPU acceleration achieves about 10x speedup on both datasets. For ADMM and fast AMA, the GPU implementation significantly outperforms the CPU implementation.

    \item \textbf{MNIST}.
    The MNIST dataset is the largest and most computationally demanding. SSNAL on GPU completes the computation in 889.87 seconds, while the CPU-based implementation requires 8928.26 seconds to solve only 95 problems, achieving more than 10x speedup.

\end{itemize}

\subsection{Performance across algorithms}
The three algorithms, SSNAL, ADMM, and fast AMA, exhibit varying levels of GPU acceleration, reflecting differences in their computational structure. 
SSNAL consistently achieves the best performance across all datasets except LUNG when implemented on GPU. On larger datasets like MNIST and COIL-20, SSNAL on GPU achieves over 10x speedups compared to CPU, demonstrating its scalability and efficiency.
ADMM also benefits significantly from GPU acceleration, especially for larger datasets. For example, on the MNIST dataset, the GPU implementation completes in 915.35 seconds, compared to 8909.42 seconds for the CPU. However, the performance gap between ADMM and SSNAL remains evident, as ADMM generally requires more iterations to converge.
Fast AMA, while not as competitive in terms of overall runtime, still demonstrates notable improvements on GPU for larger datasets. For example, on the MNIST dataset, Fast AMA on GPU completes in 3278.37 seconds, whereas the CPU implementation fails to solve any problems within the runtime threshold \( 10T \).

\begin{figure}[htbp]
    \centering

        \begin{subfigure}{0.45\textwidth}
            \includegraphics[width=\linewidth]{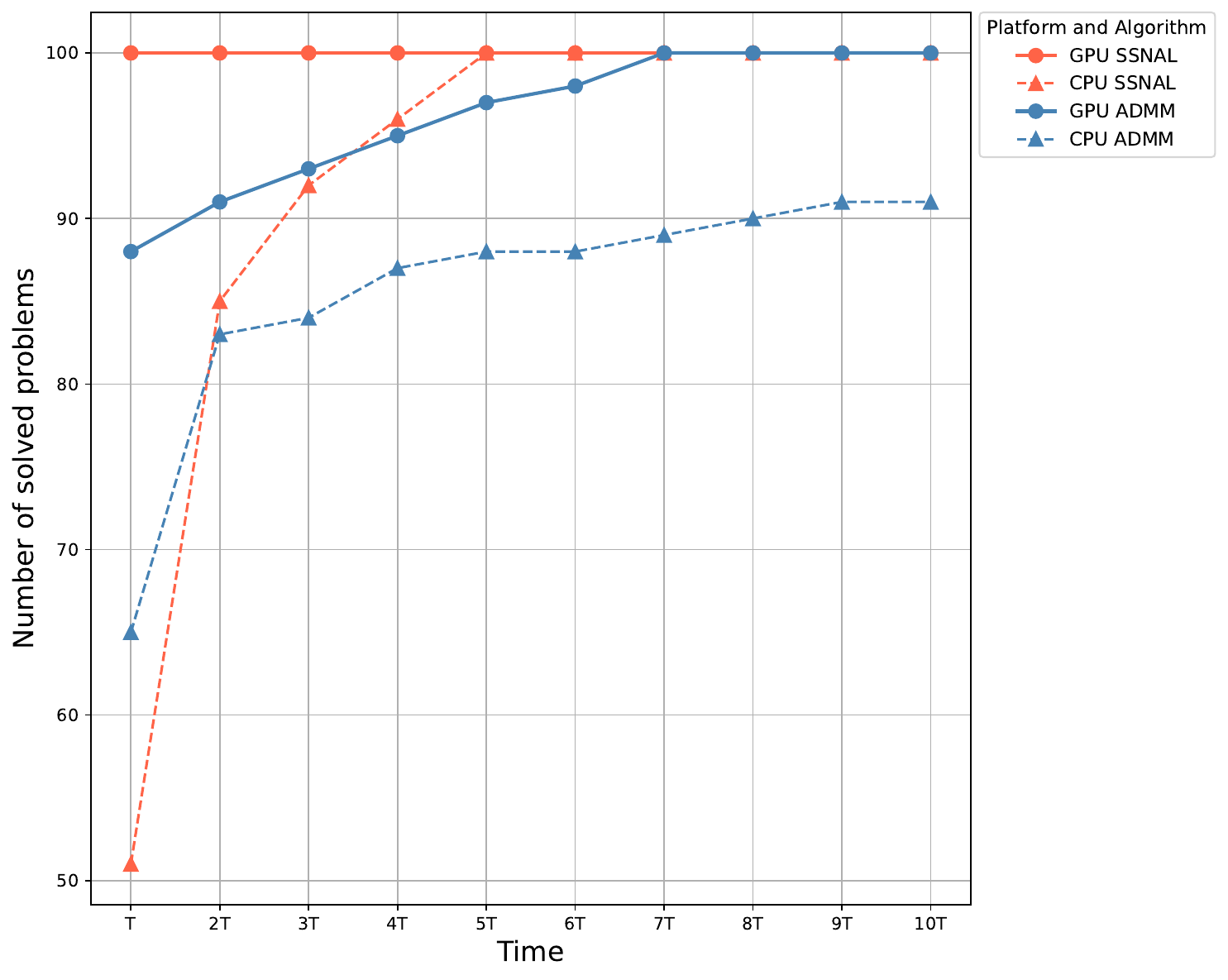}
            \caption{LIBRAS-6}
            \label{fig:LIBRAS6}
        \end{subfigure}\hspace{0.5em}
        \begin{subfigure}{0.45\textwidth}
            \includegraphics[width=\linewidth]{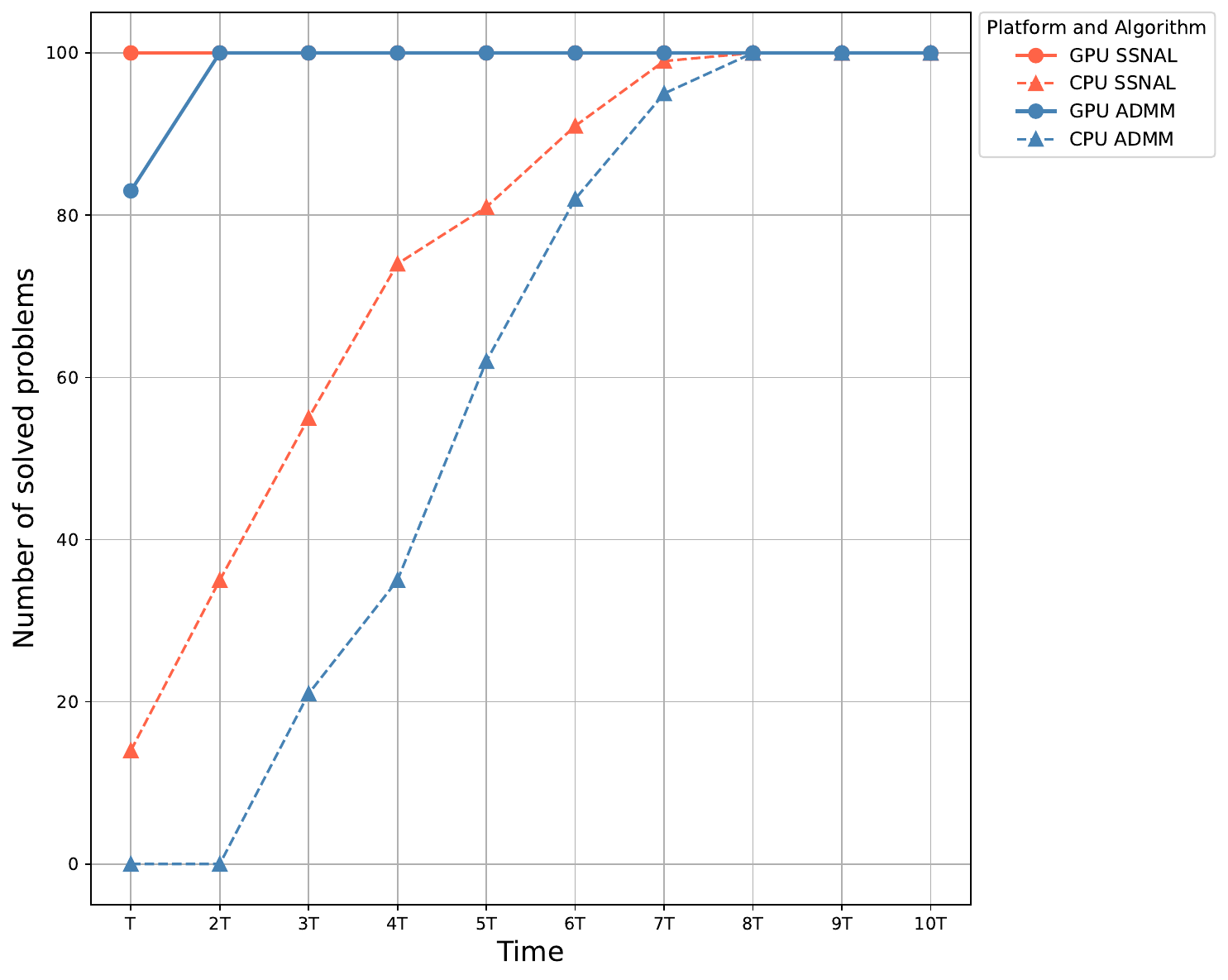}
            \caption{LIBRAS}
            \label{fig:LIBRAS}
        \end{subfigure}
    
        \vspace{0.5em} 
    
        \begin{subfigure}{0.45\textwidth}
            \includegraphics[width=\linewidth]{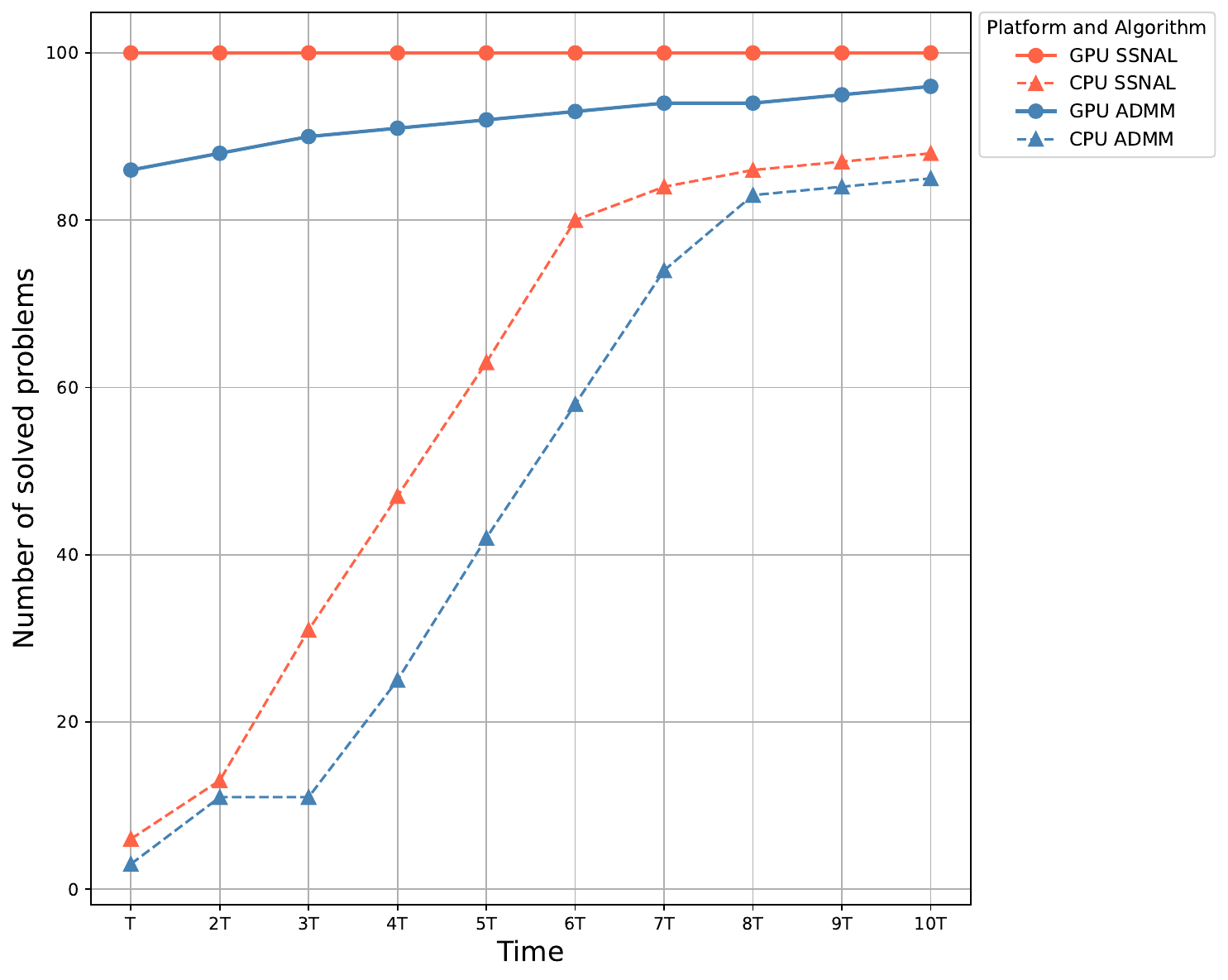}
            \caption{LUNG}
            \label{fig:lung}
        \end{subfigure}\hspace{0.5em} 
        \begin{subfigure}{0.45\textwidth}
            \includegraphics[width=\linewidth]{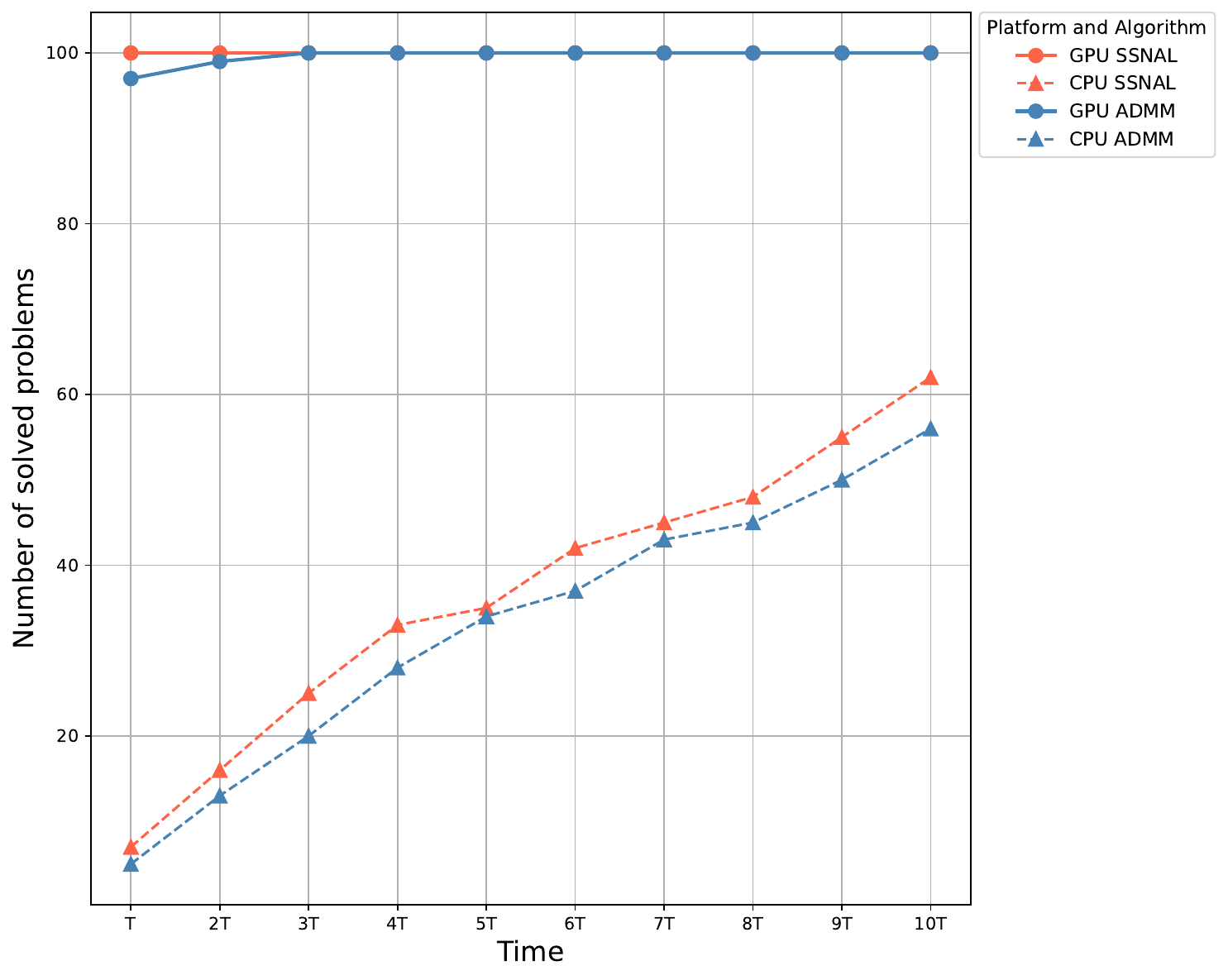}
            \caption{COIL-20}
            \label{fig:coil_profile}
        \end{subfigure}

        \vspace{0.5em}

        \begin{subfigure}{0.45\textwidth}
            \includegraphics[width=\linewidth]{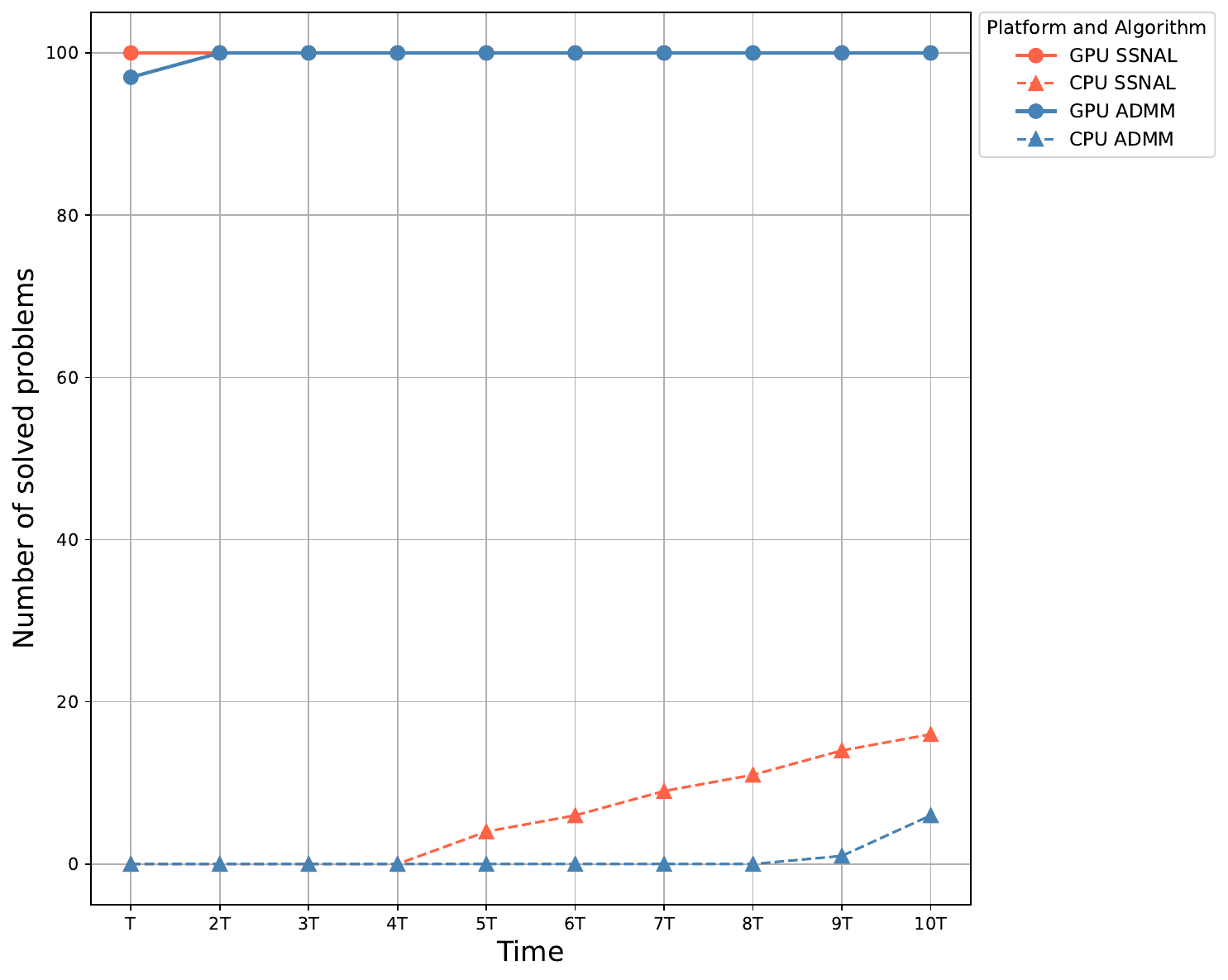}
            \caption{MNIST}
            \label{fig:MNIST}
        \end{subfigure}

    \caption{Performance profiles of different algorithms across the five benchmark datasets. The x-axis represents the normalized runtime, ranging from the baseline \( T \) to \( 10T \). The y-axis indicates the number of problems solved within the corresponding runtime.}  
    \label{fig:dataset_comparison}
\end{figure}

To provide a clearer comparison of computational efficiency across different algorithms and platforms, we plot the performance profiles of the SSNAL and ADMM on {\it PyClustrPath}(GPU) and {\it ConvexClustering}(CPU) based on the experimental results in Table~\ref{tab:computation_time_1e-6}, which are shown in Figure~\ref{fig:dataset_comparison}. The numerical results demonstrate that GPU acceleration significantly enhances the efficiency of convex clustering algorithms, particularly for larger datasets. Among the tested algorithms, the SSNAL algorithm with GPU acceleration achieves the best numerical performance on the selected datasets.   

\section{Conclusion}
This paper introduced an efficient Python package {\it PyClustrPath} for solving the convex clustering model, which implements three popular optimization algorithms and supports computation on both GPU and CPU.  By leveraging GPU acceleration, {\it PyClustrPath} significantly improves the computational efficiency, making it a new tool for generating clustering paths for large-scale datasets. There are several potential directions for us to improve {\it PyClustrPath} in the future. First, we will further optimize the implementations of the algorithms in the {\it PyClusterPath} package. We also plan to design and incorporate more efficient algorithms for solving the convex clustering model on huge-scale datasets. One of the potential algorithms is the (semi-proximal) Halpern-Peaceman-Rachford algorithm, which has achieved success in solving huge-scale linear programming problems \citep{zhang2022efficient,sun2024accelerating,zhang2024hot,chen2024hpr}. Second, we plan to implement the dimension reduction techniques introduced in \citep{yuan2021dimension,wang2023randomly} to further improve the efficiency for generating clustering paths based on the convex clustering model. Third, we plan to generalize the {\it PyClustrPath} package to other related clustering models, such as sparse convex clustering model \citep{wang2018sparse} and the convex hierarchical clustering model for graph structured data \citep{donnat2019convex}. 

\newpage
\bibliographystyle{plainnat}
\bibliography{Pyclustrpath}

\end{document}